\documentclass[10pt]{article}
\usepackage{amssymb}
\usepackage{}
\usepackage{amsmath}
\usepackage{amsfonts}
\usepackage{latexsym,amsfonts,amsmath,amsthm,makeidx,amssymb}
\usepackage{mathrsfs}
\usepackage{makeidx}
\makeindex
\newtheorem{definition}{Definition}[section]
\newtheorem{theorem}{Theorem}[section]
\newtheorem{lemma}{Lemma}[section]
\newtheorem{corollary}{Corollary}[section]
\newtheorem{proposition}{Proposition}[section]
\newtheorem{remark}{Remark}[section]

\newcommand{\R}{\mathbb R}

\newcommand{\bpp}{\begin{proposition}}
\newcommand{\epp}{\end{proposition}}

\newcommand{\bt}{\begin{theorem}}
\newcommand{\et}{\end{theorem}}
\newcommand{\bl}{\begin{lemma}}
\newcommand{\el}{\end{lemma}}
\newcommand{\bd}{\begin{definition}}
\newcommand{\ed}{\end{definition}}
\newcommand{\bc}{\begin{corollary}}
\newcommand{\ec}{\end{corollary}}
\newcommand{\bp}{\begin{proof}}
\newcommand{\ep}{\end{proof}}
\newcommand{\bx}{\begin{example}}
\newcommand{\ex}{\end{example}}
\newcommand{\bi}{\begin{exercise}}
\newcommand{\ei}{\end{exercise}}
\newcommand{\bo}{\begin{prop}}
\newcommand{\eo}{\end{prop}}
\newcommand{\br}{\begin{remark}}
\newcommand{\er}{\end{remark}}
\newcommand{\be}{\begin{equation}}
\newcommand{\ee}{\end{equation}}
\newcommand{\ba}{\begin{align}}
\newcommand{\ea}{\end{align}}
\newcommand{\bn}{\begin{enumerate}}
\newcommand{\en}{\end{enumerate}}
\newcommand{\bg}{\begin{align*}}
\newcommand{\bcs}{\begin{cases}}
\newcommand{\ecs}{\end{cases}}

\newcommand{\bean}{\begin{eqnarray*}}
\newcommand{\eean}{\end{eqnarray*}}

\allowdisplaybreaks

\numberwithin{equation}{section}

\begin{document}
\title{\bf Existence and Uniqueness of Normalized Solutions for  the Kirchhoff equation}
\date{}
\author{
{\bf Xiaoyu Zeng\quad Yimin Zhang }\thanks{ Corresponding author.
E-mail: xyzeng@whut.edu.cn (X. Y. Zeng);  zhangyimin@whut.edu.cn (Y. M. Zhang).
} \\
{\small\it  Department of Mathematics, School of Sciences, Wuhan University of Technology,
 }\\
{\small\it Wuhan 430070, P. R. China}
}

\maketitle



\vskip0.20in

\begin{abstract}
 For a class  of Kirchhoff functional, we first give a complete classification with respect to the exponent $p$  for its $L^2$-normalized critical points, and  show that the minimizer of the functional, if  exists,  is unique up to translations. Secondly, we search for the mountain pass type critical point for the functional  on the $L^2$-normalized manifold, and also prove that this type critical  point  is unique up to translations. Our proof relies only on  some simple energy estimates and  avoids using the concentration-compactness principles. These conclusions  extend some known results in previous papers.
 \end{abstract}



\vskip0.2in


\section{Introduction}
\indent\indent In this paper, we study the following   Kirchhoff equation
\begin{equation}\label{eq1.1}
-\Big(a+b\int_{\mathbb{R}^N}|\nabla u|^2dx\Big)\Delta u-|u|^pu= \mu u\text{ in }\ \R^N,
\end{equation}
where $a,b>0, 1 \le N\le3$, and $0<p<2^*-2$ with $2^*=+\infty$ if $N=1,2$, or $2^*=6$ if $N=3$.
Equation (\ref{eq1.1}) is related to the stationary solutions of
\begin{equation}\label{eeq1.2}
u_{tt}-\Big(a+b\int_{\mathbb{R}^N}|\nabla u|^2dx\Big)\Delta u=f(x,u),
\end{equation}
where $f(x,u)$ is a general nonlinearity.  The problem (\ref{eeq1.2}) was proposed by Kirchhoff \cite{Ki} and models free vibrations of elastic
strings by  taking  into account the changes in length of the string produced by transverse vibrations. Comparing  with the semilinear equations (i.e., setting  $b=0$ in above two equations), it is much more challenge and interesting   to investigate equations (\ref{eq1.1}) and (\ref{eeq1.2}) in view of the existence of the nonlocal term $\big(\int_{\R^N}|\nabla u|^2dx\big) \triangle u$.

After the pioneering work of \cite{Po} and \cite{L}, much attention was paid to these above two equations. For instance, replacing the term $|u|^pu$ with a general nonlinearity $f(x,u)$, there are many results on the existence of solutions for equation (\ref{eq1.1}), one can refer \cite{DPS,HLP,HZ} and the references therein. Equation  (\ref{eq1.1}) can be viewed as an eigenvalue problem by taking $\mu$ as an unknown Lagrange multiplier. From this point of view, one can solve (\ref{eq1.1}) by studying some constrained variational problems and obtain normalized solutions. Motivated by the works of \cite{BJL,Y1}, we consider the following minimization problem:
\begin{equation}\label{eq1.2}
I(c):=\inf_{u\in S_c}E(u)
\end{equation}
where
\begin{equation}\label{eq1.0}
E(u)=\frac{a}{2}\int_{\mathbb{R}^N}|\nabla u|^2dx+\frac{b}{4}\big(\int_{\mathbb{R}^N}|\nabla u|^2dx\big)^2-\frac{1}{p+2}\int_{\mathbb{R}^N}|u|^{p+2}dx,
\end{equation}
and
$$S_c:=\big\{u\in H^1(\R^N): \int_{\mathbb{R}^N} |u|^2dx=c^2\big\}.$$
\begin{remark}
If $u\in S_c$ is a minimizer of problem (\ref{eq1.2}), then there exists $\mu\in\R$ such that $E'(u)=\mu u$, namely, $u\in S_c$ is a solution of (\ref{eq1.1}) for some $\mu\in \R$. Therefore, it is natural  to obtain normalized solutions  by investigating problem (\ref{eq1.2}).

\end{remark}

 Ye in \cite{Y1} proved that when $0<p<\frac{8}{N}$, there exists $c_p^*>0$ such that  (\ref{eq1.2}) has minimizers if and only if $c>c_p^*$ and $0<p\leq \frac{4}{N}$, or $c\geq c_p^*$  and $\frac{4}{N}<p<\frac{8}{N}$,  where $c_p^*$ is given by
 \begin{equation}c^*_p=\begin{cases}
 0, \ \ &0<p<\frac{4}{N};\\
 a^\frac{4}{N}|Q|_{L^2}, &p=\frac{4}{N};\\
 \inf\{c\in(0,+\infty): I(c)<0\}, &\frac{4}{N}<p<\frac{8}{N};
 \end{cases}
\end{equation}
and  $Q(|x|)$ is the unique (up to translations) radially symmetric positive solution of the following equation  in $%
H^{1}(\mathbb{R}^{N})$ :
\begin{equation}\label{1.3}
-\Delta u+\frac{4+2p-Np}{Np}u-\frac{4}{Np}|u|^{p}u=0,\ 0<p< 2^*-2.
\end{equation}%
While for the case of $\frac{8}{N}\leq p<2^*-2$, problem (\ref{eq1.2}) cannot be attained.
Recently, Guo and Wang in \cite{GW} gave the explicit form of $c^*_p$ for $p\in(\frac{4}{N},\frac{8}{N})$. We note that the arguments of \cite{Y1} mainly depend on the application of the concentration-compactness principle. By ruling out the cases of {\em vanishing} and {\em dichotomy }, the author obtained the compactness of a minimizing sequence.

 In what follows, by observing the special form of the functional (\ref{eq1.0}), we intend to give a new proof for the above results  for problem (\ref{eq1.2}) in a simple way, where only  technical energy estimates are involved and the concentration-compactness principle is avoided. Especially, our arguments also show that the minimizer of $(\ref{eq1.0})$, if  exists, is unique and must be a scaling of $Q(x)$.
Before stating our main results, we first recall that, combining with the  Pohozaev and Nehari identity ,  $Q(x)$ satisfies
\begin{equation}  \label{eq2.1a}
\int_{\mathbb{R}^N}|\nabla Q(x)|^2dx=\int_{\mathbb{R}^N}|
Q(x)|^2dx=\frac{2}{p+2}\int_{\mathbb{R}^N}| Q(x)|^{p+2}dx.
\end{equation}
Moreover, $Q$ is an optimizer of the following sharp Gagliardo-Nirenberg
inequality \cite{W}
\begin{equation}  \label{1.5}
\int_{\mathbb{R}^N}|u|^{p+2}dx\leq \frac{p+2}{2|Q|_{L^2}^p}\Big(
\int_{\mathbb{R}^N}|\nabla u|^2dx\Big)^{\frac{Np}{4}}\Big(\int_{\mathbb{R}^N}|u|^2dx\Big)^{\frac{2(p+2)-Np}{4}}, \forall \ u\in H^1(\mathbb{R}^3).
\end{equation}
We remark that, similar to \cite{W}, one can prove that  all optimizers of (\ref{1.5}) are indeed  the scalings and translations  of $Q(x)$, i.e., belong to the following set
\begin{equation}\label{eq1.7}
\{\lambda Q(\alpha x+y):\ \alpha, \beta\in \R^+, y\in \R^N\}.
\end{equation}

Our first  theorem gives the existence and uniqueness of minimizers for problem (\ref{eq1.2}).
\begin{theorem}\label{thm1}
\begin{itemize}
  \item [(i)] When $0<p<\frac{4}{N}$,  problem (\ref{eq1.2}) has a unique minimizer (up to translations), which is the form of
  $u_c=\frac{c\lambda_p^\frac{N}{2}}{|Q|_{L^2}}Q(\lambda_p x),$
  where $\lambda_p=\frac{t_p^\frac{1}{2}}{c}$ with $t_p$ being the unique minimum point of the function
  \begin{equation}\label{eq1.8}
  f_p(t)=\frac{a}{2}t+\frac{b}{4}t^2-\frac{c^\frac{2(p+2)-Np}{2}}{2|Q|_{L^2}^p}t^\frac{Np}{4},\ t\in(0,+\infty).
  \end{equation}
  \item [(ii)] When $p=\frac{4}{N}$,  problem (\ref{eq1.2}) has no minimizer if $c\leq a^\frac{N}{4}|Q|_{L^2}$. On the contrary, if $c> a^\frac{N}{4}|Q|_{L^2}$, then (\ref{eq1.2}) has a unique  minimizer (up to translations)
      \begin{equation}\label{eq1.9}
      u_c=\frac{c\lambda_p^\frac{N}{2}}{|Q|_{L^2}}Q(\lambda_p x) \text{ where }\lambda_p=\frac{(c^\frac{4}{N}-a|Q|_{L^2}^\frac{4}{N})^\frac{1}{2}}{c|Q|_{L^2}^\frac{2}{N}b^\frac{1}{2}}.\end{equation}
      Also,$
      I(c)=-\frac{(c^\frac{4}{N}-a|Q|_{L^2}^\frac{4}{N})^2}{4b|Q|_{L^2}^\frac{8}{N}}.$
  \item [(iii)] When $\frac{4}{N}<p<\frac{8}{N}$, problem (\ref{eq1.2}) has no minimizer if
  \begin{equation}\label{eq1.11}
  c<c^*:=\Big(2|Q|_{L^2}^p\big(\frac{2a}{8-Np}\big)^\frac{8-Np}{4}\big(\frac{b}{Np-4}\big)^\frac{Np-4}{4}\Big)^\frac{2}{2(p+2)-Np}.\end{equation}
  On the other hand, if $c\geq c^*$, then (\ref{eq1.2}) has a unique minimizer (up to translations)
  \begin{equation}\label{eq1.12}
  u_c=\frac{c\lambda_p^\frac{N}{2}}{|Q|_{L^2}}Q(\lambda_p x) \text{ with }\lambda_p=\frac{1}{c}\Big(\frac{2(Np-4)a}{(8-Np)b}\Big)^\frac{1}{2}.\end{equation}
  Moreover, we have
$
  I(c)=\frac{{c^*}^\frac{2(p+2)-Np}{2}-c^\frac{2(p+2)-Np}{2}}{2|Q|_{L^2}^p}\Big(\frac{2(Np-4)a}{(8-Np)b}\Big)^\frac{Np}{4}\ \text{ for any }c\geq c^*.
$

  \item [(iv)] When $p\geq\frac{8}{N}$, problem (\ref{eq1.2}) has no minimizer for all $c>0$.

\end{itemize}
\end{theorem}
Theorem \ref{thm1} tells that the minimizer of (\ref{eq1.2}) must be a scaling of $Q(x)$, which extends \cite[Theorem 1.1]{Y1}, where the existence of minimizers for (\ref{eq1.2}) was discussed. From the above theorem, we see that problem (\ref{eq1.2}) has no minimizer if $p\geq\frac{8}{N}$. Thus, to obtain the normalized solutions for (\ref{eq1.1}), one may   search for saddle point for functional (\ref{eq1.0}). Stimulated by \cite{J1,BJL}, we investigate the  mountain pass type critical point  for  $E(\cdot)$  on $S_c$.
\begin{definition}
Given $c>0$, the functional $E(\cdot)$ is said to have mountain pass geometry on $S_c$ if there exists $K(c)>0$ such that \begin{equation}\label{eq1.14}\gamma(c):=\inf_{h\in\Gamma(c)}\max_{t\in[0,1]}E\big(h(t)\big)>\max\{E(h(0)), E(h(1))\}
\end{equation}
holds in the set
$
\Gamma(c)=\Big\{h\in C\big([0,1]; S_c\big)| h(0)\in A_{K(c)} \text{ and } E\big(h(1)\big)<0\Big\},
$
where $A_{K(c)}=\{u\in S_c: |\nabla u|_{L^2}^2\leq K(c)\}$.
\end{definition}

By  studying some analytic properties of $\gamma(c)$ and involving rigorous arguments , Ye in \cite{Y1,Y2} proved respectively that,

\begin{equation}\label{eq1.18}
\text{ either } p>\frac{8}{N}, \text{ or } p=\frac{8}{N} \text{ and } c>c_*:=\big(\frac{b|Q|_{L^2}^\frac{8}{N}}{2}\big)^\frac{N}{8-2N},
 \end{equation}
$E(\cdot)$ possesses the mountain pass geometry on $S_c$. Moreover, there exists $u_c\in S_c$ such that $E(u_c)=\gamma(c)$, and $u_c$ is a solution of (\ref{eq1.1}) with some $\mu\in \R^-$. Immediately, for the second case, Ye in \cite{Y3} further  studied the asymptotic behavior of $u_c$ as $c\to (c_*)^+$. Motivated by the these results and the proof of our first theorem, we attempt to investigate some  properties of problem (\ref{eq1.14}) by  introducing some new observations and  energy estimates. Moreover, as a byproduct, we show that  if  $u_c\in S_c$ is critical point of $E(\cdot)$ on the level  $\gamma(c)$, then  it is unique and indeed  a scaling of $Q(x)$.
 Still let $f_p(\cdot)$ be given by (\ref{eq1.8}) and note  that it has a  unique maximum point in $(0,+\infty)$ once (\ref{eq1.18}) is assumed. Then, we have the following theorem.

\begin{theorem}\label{thm2}
Assume (\ref{eq1.18}) holds and let  $\bar t_p$ be the  unique maximum point of $f_p(t)$ in $(0,+\infty)$. Then
$
\gamma(c)=f_p(\bar t_p)
$
and it can be attained by
$
\bar u_c=\frac{c\bar\lambda_p^\frac N2}{|Q|_{L^2}}Q(\bar \lambda_p x),
$
where $\bar \lambda_p=(\bar t_p)^\frac12/c$. Also, $\bar u_c$ is a solution of (\ref{eq1.1}) for some $\mu\in\R^-$. Moreover, $\bar u_c$ is the unique solution of (\ref{eq1.14}) in the following sense: if $\bar u\in S_c$ is a critical point of $E(\cdot)$ on $S_c$ and its energy equals to $\gamma(c)$, namely,
\begin{equation}\label{eq1.19}
E'(\bar u)|_{S_c}=0\ \text{ and }E(\bar u)=\gamma(c).\end{equation}
Then, up to translations, $\bar u=\bar u_c$.
\end{theorem}
\begin{remark}
If $p=\frac{8}{N}$ and $c>c_*$, one can easily check that $\bar t_p=\frac{a}{b}[\big(\frac{c}{c_*}\big)^\frac{8-2N}{N}-1]^{-1}$,
we thus deduce from Theorem \ref{thm2} and (\ref{eq1.18}) that
$\bar u_c
=\Big(\frac{a^2(cc_*)^\frac{8-4N}{N}}{2bc_*^2(\big(\frac{c}{c_*}\big)^\frac{8-2N}{N}-1)^2}\Big)^\frac{N}{8}Q\Big(\big[\frac{a}{bc^2(\big(\frac{c}{c_*}\big)^\frac{8-2N}{N}-1)}\big]^\frac{1}{2}x\Big)
$
and
$\gamma(c)=\frac{1}{4b}[\big(\frac{c}{c_*}\big)^\frac{8-2N}{N}-1]^{-1}.$
This extends the results of \cite{Y3} where the asymptotic behaviors  of $\bar u_c$ and the value of $\gamma (c)$ as $c\to (c_*)^+$ were studied.
\end{remark}

\section{Proof of  Main Results. }
In this section, we give the proof of Theorems \ref{thm1} and \ref{thm2} by employing the Gagliardo-Nirenberg
inequality (\ref{1.5}) and some energy estimates. We first remark that by a simple rescaling, one can easily show that
\begin{equation}\label{eq2.1}
I(c)\leq0\ \text{ for all } c>0 \text{ and } 0<p<2^*-2.
\end{equation}
Moreover, utilizing (\ref{1.5}), we see that for any $u\in S_c$,
\begin{align}
E(u)&\geq\frac{a}{2}\int_{\mathbb{R}^N}|\nabla u|^2dx+\frac{b}{4}\Big(\int_{\mathbb{R}^N}|\nabla u|^2dx\Big)^2-\frac{c^\frac{2(p+2)-Np}{2}}{2|Q|_{L^2}^p}\Big(
\int_{\mathbb{R}^N}|\nabla u|^2dx\Big)^{\frac{Np}{4}}\nonumber\\
&=f_p(t)\ \ \text{ by setting }t=\int_{\mathbb{R}^N}|\nabla u|^2dx.\label{eq2.2}
\end{align}
where $f_p(\cdot)$ is given by (\ref{eq1.8}).

\noindent\textbf{\noindent{Proof of Theorem \ref{thm1}}.}
\textbf{(i).}Since $p<\frac{4}{N}$, one can easily check that $f_p(t)\ (t\in(0,\infty))$ attains its minimum at a unique point, denoted by $t_p$. Therefore, we obtain from (\ref{eq2.2}) that
\begin{equation}\label{eq2.3}
I(c)=\inf_{u\in S_c}E(u)\geq f_p(t_p).\end{equation}
On the other hand, set
\begin{equation}\label{eq2.4}
u_\lambda(x)=\frac{c\lambda^\frac{N}{2}}{|Q|_{L^2}}Q(\lambda x),\end{equation}  where $\lambda>0 $ will be determined later.
Then, $u_\lambda\in S_c$ and it follows from  (\ref{eq2.1a}) that
$$\int_{\mathbb{R}^N} |\nabla u_\lambda|^2dx=c^2\lambda^2; \ \ \int_{\mathbb{R}^N}|u_\lambda|^{p+2}dx=\frac{(p+2)c^{p+2}\lambda^\frac{Np}{2}}{2|Q|_{L^2}^p}.$$
Consequently,
\begin{equation}\label{eq2.44}
\begin{split}
E(u_\lambda)&=\frac{a}{2}c^2\lambda^2+\frac{b}{4}(c^2\lambda^2)^2-\frac{c^\frac{2(p+2)-Np}{2}}{2|Q|_{L^2}^p}(c^2\lambda^2)^\frac{Np}{4}=f_p(c^2\lambda^2).
\end{split}
\end{equation}
Choosing $\lambda=t_p^\frac{1}{2}/c$, i.e., $c^2\lambda^2=t_p$, it follows from  (\ref{eq2.44}) that
 $I(c)\leq E(u_\lambda)=f_p(t_p).$
Together with (\ref{eq2.3}), we deduce that
\begin{equation}\label{eq2.5}
I(c)=f_p(t_p)=\inf_{t\in\R^+}f_p(t),\end{equation}
 and
$u_\lambda$ with  $\lambda=t_p^\frac{1}{2}/c$, i.e., $u_\lambda =u_c=\frac{c}{|Q|_{L^2}}{\left(t_p^\frac{1}{2}/c\right)}^\frac{N}{2}Q(t_p^\frac{1}{2}x/c )$ is a minimizer of (\ref{eq1.2}).

It remains to prove that, up to translations, $u_c$ is the unique minimizer of (\ref{eq1.2}). Indeed, if $u_0\in S_c$ is a minimizer, it then follows from (\ref{eq2.2}) that
$
E(u_0)\geq f_p(t_0), \ \text{ with } t_0:=\int_{\mathbb{R}^N}|\nabla u_0|^2dx,
$
where the ``$=$" holds if and only if $u_0$ is an optimizer of (\ref{1.5}).
This and  (\ref{eq2.5}) further imply that $t_0=t_p$ and
$f_p(t_0)=E(u_0)$. Thus,  $u_0$ is an optimizer of (\ref{1.5}) and it  follows from (\ref{eq1.7}) that up to translations, $u_0$ must be the form of
$u_0(x)=\alpha Q(\beta x).$
Utilizing $\int_{\mathbb{R}^N}|u_0|^2dx=c^2$, $\int_{\mathbb{R}^N}|\nabla u_0|^2dx=t_p$ and (\ref{eq2.1a}), we see that
$\alpha=\frac{c}{|Q|_{L^2}}\big(\frac{t_p^\frac{1}{2}}{c}\big)^\frac{N}{2} \text{ and } \beta=\frac{t_p^\frac{1}{2}}{c},$
hence, $u_0=u_c$.

\textbf{(ii).} Since $p=\frac{4}{N}$, then,
\begin{equation}\label{eq2.7}
f_p(t)=\frac{1}{2}\Big(a-\frac{c^\frac{4}{N}}{|Q|_{L^2}^\frac{4}{N}}\Big)t+\frac{b}{4}t^2.\end{equation}
 If $c\leq a^\frac{N}{4}|Q|_{L^2}$ one can easily deduce from (\ref{eq2.2}) that
$
E(u)\geq f\big(\int_{\mathbb{R}^N}|\nabla u|^2dx\big)>0 \ \text{ for all } u\in S_c.
$
In view of (\ref{eq2.1}), this indicates that (\ref{eq1.2}) has no minimizer.
Next, we turn to the case of $c>a^\frac{N}{4}|Q|_{L^2}$. From (\ref{eq2.7}), we know  that  $f_p(t)(t\in(0,+\infty))$ attains its minimum at the unique point
$t_p=\frac{c^\frac{4}{N}-a|Q|_{L^2}^\frac{4}{N}}{b|Q|_{L^2}^\frac{4}{N}b}.$
Similar to the arguments of part (i), one can prove that, up to translations $u_c$ given by (\ref{eq1.9}) is the unique minimizer of (\ref{eq1.2}) and  the energy  $I(c)=f_p(t_p)=-\frac{(c^\frac{4}{N}-a|Q|_{L^2}^\frac{4}{N})^2}{4b|Q|_{L^2}^\frac{8}{N}}$ .

\textbf{(iii).} For the case $\frac{4}{N}<p<\frac{8}{N}$, let
$\alpha=\frac{8-Np}{4} \ \text{and }\beta=1-\alpha=\frac{Np-4}{4},$
it follows from the Young's inequality that, for any $t>0$,
\begin{equation*}
\begin{split}
\frac{a}{2}t+\frac{b}{4}t^2=\alpha\left(\frac{a}{2\alpha}t\right)+\beta\left(\frac{b}{4\beta}t^2\right)\geq \left(\frac{a}{2\alpha}\right)^\alpha\left(\frac{b}{4\beta}\right)^\beta t^{\alpha+2\beta}=\big(\frac{2a}{8-Np}\big)^\frac{8-Np}{4}\big(\frac{b}{Np-4}\big)^\frac{Np-4}{4} t^\frac{Np}{4}
\end{split}
\end{equation*}
where the ``$=$" in the second inequality holds if and only if
$
\frac{a}{2\alpha}t=\frac{b}{4\beta}t^2, \ \text{ i. e., } t=t_0:=\frac{2\beta a}{\alpha b}=\frac{2(Np-4)a}{(8-Np)b}.
$
In view of (\ref{eq2.2}) and noting that $c_*$ is given by  (\ref{eq1.11}), we therefore have
\begin{equation}\label{eq2.10}
E(u)\geq \frac{{c^*}^\frac{2(p+2)-Np}{2}-c^\frac{2(p+2)-Np}{2}}{2|Q|_{L^2}^p}t_0^\frac{Np}{4}=f_p(t_0) \text{ for all } u\in S_c.
\end{equation}

If $c<c^*$, we then deduce from (\ref{eq2.10}) that $E(u)>0$ for all $u\in S_c$. Thus, problem (\ref{eq1.2}) cannot be achieved
for (\ref{eq2.1}).

If $c\geq c^*$, on the one hand, we deduce from (\ref{eq2.10})
that
$
I(c)\geq f_p(t_0)$.
On the other hand, let $u_\lambda(x)$ be as in (\ref{eq2.4}) and set $\lambda=t_0^\frac{1}{2}/c$, then
$I(c)\leq E(u_\lambda)=f_p(t_0)$.
This indicates that $u_\lambda$ is a minimizer of (\ref{eq1.2}) and $I(c)=f_p(t_0)=\frac{{c^*}^\frac{2(p+2)-Np}{2}-c^\frac{2(p+2)-Np}{2}}{2|Q|_{L^2}^p}\Big(\frac{2(Np-4)a}{(8-Np)b}\Big)^\frac{Np}{4}$ for any $c\geq c^*$.
The uniqueness of minimizers can be proved by the same argument of part (i).

\textbf{(iv).} If $p>\frac{8}{N}$, or $p=\frac{8}{N}$ and $c>\big(\frac{b|Q|_{L^2}^\frac{8}{N}}{2}\big)^\frac{N}{8-2N}$, it follows (\ref{eq2.4}) and (\ref{eq2.44}) that
$I(c)\leq\lim_{\lambda\to+\infty}E(u_\lambda)=-\infty,$
and thus problem (\ref{eq1.2}) cannot be attained. On the other hand, if $p=\frac{8}{N}$ and $c\leq\big(\frac{b|Q|_{L^2}^\frac{8}{N}}{2}\big)^\frac{N}{8-2N}$,  from  (\ref{eq2.2}) we have
$E(u)>0 \ \text{ for all } u\in S_c$.
This together with (\ref{eq2.1}) obviously indicates that problem (\ref{eq1.2}) cannot be attained.
\qed

\vskip0.1truein

\noindent\textbf{Proof of Theorem \ref{thm2}.}
Firstly, similar to the proof of \cite[Lemma 3.1]{Y2}, from the definition 1.1, one can prove that there exists $K(c)>0$  which  can be chosen small enough such that,  $E(\cdot)$ admits mountain pass geometry on $S_c$ if (\ref{eq1.18}) is assumed. In what follows, we thus always assume that $K(c)<\bar t_p$, where $\bar t_p$ denotes  the unique maximum point of $f_p(t)$ in $(0,+\infty)$.

For any $h(s)\in \Gamma(c)$, one can deduce from  (\ref{eq2.2}) that
\begin{equation}\label{eq2.13}
E(h(s))\geq f_p\Big(\int_{\mathbb{R}^N}|\nabla h(s)|^2dx\Big),
\end{equation}
where ``$=$" holds if and only if $h(s)\in S_c$ is an optimizer of (\ref{1.5}), i. e., up to translations,
\begin{equation}\label{eq2.11}
(h(s))(x)=\frac{c\alpha^\frac{N}{2}}{|Q|_{L^2}}Q(\alpha x)\text{ for some }\alpha>0.\end{equation}
Since $h(0)\in A_{K(c)}$ with $K(c)<\bar t_p$, and note that  $f_p(t)>0\ \forall t\in(0,\bar t_p]$, we thus have
\begin{equation}\label{eq2.155}
\int_{\mathbb{R}^N}|\nabla h(0)|^2dx<\bar t_p<\int_{\mathbb{R}^N}|\nabla h(1)|^2dx.
\end{equation}
As a consequence of (\ref{eq2.13}) and (\ref{eq2.155}), there holds that
\begin{equation}\label{eq2.15}\max_{s\in[0,1]}E\big(h(s)\big)\geq f_p(\bar t_p)=\max_{t\in\R^+}f_p(t).\end{equation}
Thus,
\begin{equation}\label{eq2.17}
\gamma(c)\geq f_p(\bar t_p).\end{equation}
On the contrary, let  $u_\lambda(x)$ be the trial function given by (\ref{eq2.4}) with $\lambda=\bar \lambda_p=(\bar t_p)^\frac12/c$.
Set $\bar h(s):=s^\frac{N}{4}u_\lambda(s^\frac{1}{2}x)$, then one can check that $E(\bar h(s))=f_p(\bar t_p s)$. Choosing $0<\tilde t_p<\bar t_p$ small enough such that $\bar h(\tilde t_p /\bar t_p)\in A_{K(c)}$, and  $\hat t_p>\bar t_p$ such that $f_p(\hat t_p)<0$, let
$h(s)=\bar h((1-s)\tilde t_p/\bar t_p+\hat t_p s/\bar t_p).$
Then,
$h(0)=\bar h(\tilde t_p /\bar t_p)\in A_{K(c)}\text{ and  }E(h(1))=E(\bar h(\hat t_p/\bar{t}_p))=f_p(\hat t_p)<0.$
This indicates that  $h\in\Gamma(c)$, and
$$\gamma(c)\leq \max_{t\in[0,1]}E\big(h(t)\big)=E(u_{\bar\lambda_p})=f_p(\bar t_p).$$
Combing with (\ref{eq2.17}), we deduce that  $\gamma(c)=f_p(\bar t_p)$ and $u_{\bar \lambda_p}=\bar u_c(x)=\frac{c}{|Q|_{L^2}}{\left(\bar t_p^\frac{1}{2}/c\right)}^\frac{N}{2}Q(\bar t_p^\frac{1}{2}x/c )$ is a solution of problem (\ref{eq1.14}).

We next prove that $\bar u_c$ satisfies equation (\ref{eq1.1}) for some $\mu\in\R^-$. Actually, in view of $f_p'(\bar t_p)=0$ and $\bar\lambda_p=(\bar t_p)^\frac12/c$, we have

\begin{equation}\label{eq2.188}\frac{Npc^\frac{2(p+2)-Np}{2}}{4|Q|_{L^2}^p}(\bar t_p)^\frac{Np-4}{4}=a+b\bar t_p=a+b\int_{\mathbb{R}^N}|\nabla \bar u_c|^2dx.\end{equation}
Moreover, since $Q(x)$ is a solution of (\ref{1.3}) and note that $\bar\lambda_p=(\bar t_p)^\frac12/c$, it follows that $\bar u_c$ satisfies
$$-\frac{Npc^\frac{2(p+2)-Np}{2}}{4|Q|_{L^2}^p}(\bar t_p)^\frac{Np-4}{4}\Delta \bar u_c-|\bar u_c|^p\bar u_c=-\frac{(4+2p-Np)(c\bar\lambda_p^\frac{N}{2})^p}{4|Q|_{L^2}^p}\bar u_c.$$
This together with (\ref{eq2.188}) indicates that $\bar u_c$ is a solution of (\ref{eq1.1}) with $\mu=-\frac{(4+2p-Np)(c\bar\lambda_p^\frac{N}{2})^p}{4|Q|_{L^2}^p}$.

We finally prove that, up to translations, $\bar u_c$ is the unique solution of $\gamma (c)$. Suppose $\bar u$ is a solution of $\gamma(c)$ and satisfies (\ref{eq1.19}), then there exists $\mu\in\R$ such that $E'(u)=\mu u$, it then follows from the Nehari and Pohozaev identity (see e.g., \cite[Lemma 2.1]{JL}) that
\begin{equation}\label{eq2.18}
\frac{a}{2}\int_{\mathbb{R}^N}|\nabla \bar u|^2dx+\frac{b}{2}\big(\int_{\mathbb{R}^N}|\nabla \bar u|^2dx\big)^2-\frac{pN}{4(p+2)}\int_{\mathbb{R}^N}|\bar u|^{p+2}dx=0.
\end{equation}
 Let $\hat h(s):=s^\frac{N}{4}\bar u(s^\frac{1}{2}x)$ and
 \begin{equation}
 g(s):=E(\hat h(s))=\frac{as}{2}\int_{\mathbb{R}^N}|\nabla \bar u|^2dx+\frac{bs^2}{4}\big(\int_{\mathbb{R}^N}|\nabla \bar u|^2dx\big)^2-\frac{s^\frac{pN}{4}}{p+2}\int_{\mathbb{R}^N}|\bar u|^{p+2}dx.
 \end{equation}
 (\ref{eq2.18}) indicates  that $g(s)$ $(s\in(0,+\infty))$ attains its maximum at the unique point $s=1$, and $\lim_{s\to+\infty}g(s)=-\infty$. Choosing $0<\tilde s<1<\hat s$  such that $\hat h(\tilde s)\in A_{K(c)}$ and $g(\hat s)<0$, then $h_0(s):=\hat h((1-s)\tilde s+\hat ss)\in \Gamma(c)$ and $\max_{s\in[0,1]}E(h_0(s))=E(\bar u)$. As the arguments of (\ref{eq2.13}) and (\ref{eq2.15}), we see that
$$f_p(\bar t_p)=\gamma(c)=E(\bar u)=\max_{s\in[0,1]}E(h_0(s))\geq \max_{t\in\R^+}f_p(t)=f_p(\bar t_p). $$
 Together with (\ref{eq2.11}), this means  that  $\bar u$ must be the form of $\frac{c\alpha^\frac{N}{2}}{|Q|_{L^2}}Q(\alpha x)$ for some $\alpha>0$. Take it into the equality $f_p(\bar t_p)=E(\bar u)$, we further obtain that $\alpha=\bar\lambda_p$ and  $\bar u=\bar u_c$.
\qed

\vskip 0.1truein

\noindent {\bf Acknowledgements:} X. Y. Zeng is  supported by NSFC grant 11501555, and  Y. M. Zhang is  supported by NSFC grant 11471330. This work is also partially supported by the Fundamental Research Funds for the Central Universities(WUT: 2017 IVA 075 and 2017 IVA 076).

\end{document}